\documentclass[a4paper,10pt]{article}
\usepackage[utf8]{inputenc}

\usepackage[all]{xy}
\usepackage{amsmath,amssymb,amsfonts,amsthm,graphicx}
\usepackage[margin=1.2in]{geometry}

\newtheorem{thm}{Theorem}[section]

\newtheorem{cor}[thm]{Corollary}
\newtheorem{prop}[thm]{Proposition}
\newtheorem{df}[thm]{Definition}
\newtheorem{lm}[thm]{Lemma}

\newcommand{\F}{\mathbb{F}}

\title{Recognising Abelian Sylow Subgroups in Finite Groups}
\author{Julian Brough}

\begin{document}
\date{}
\maketitle

\begin{center}
\small
\textit{Faculty of Mathematics, Centre for Mathematical Sciences,}

\textit{Wilberforce Road, Cambridge, England, CB3 0WA}
\end{center}
\normalsize
\begin{abstract}
Let $p$ be a prime. We prove that if a finite group $G$ has non-abelian Sylow $p$-subgroups, and the class size of every $p$-element in $G$ is coprime to $p$; then $G$ contains a simple group as a subquotient which exhibits the same property. 
In addition we provide a list of all the simple groups and primes such that the Sylow $p$-subgroups are non-abelian and all $p$-elements have class size coprime to $p$.
This provides a solution to the problem which was left remaining after Tiep and Navarro established that for $p\not= 3,5$ this can never happen \cite{Tiep}.  
\end{abstract}

\section{Introduction}

If we are given the character table for a finite group $G$, a natural question to ask is; how much information can we obtain about the group from the table? 
In particular Richard Brauer asked if it was possible to determine if the group has abelian Sylow $p$-subgroups \cite[Problem 12]{Brauer}. 
This question was shown to be true for all primes $p$ in \cite{Sandling}.
In addition to this result, Camina and Herzog established that for the prime $p=2$, we do not require the whole character table to determine if the Sylow 2-subgroup is abelian; in fact it is enough to use the conjugacy class sizes and show all $2$-elements have class size coprime to $2$ \cite{CamHer}. 
This prompts us to ask the question, is it possible to tell if a Sylow $p$-subgroup is abelian if we know all the class sizes of $p$-elements are coprime to $p$?

We give a solution to this problem by studing groups with the following property.
\begin{df}
 Let $G$ be a finite group, and $p$ a prime. Then we call $G$ a $cl_p$ group if the Sylow $p$-subgroups of $G$ are non-abelian and the class size of every $p$-element in $G$ is not divisible by $p$.
\end{df}

\begin{thm}[\textbf{Main Theorem}]\label{MainTheorem}
Let $G$ be a finite group, and $p$ a prime. If $G$ is a $cl_p$ group, then $p=3$ or $5$ and $G$ must contain one of the following simple $cl_p$ groups as a subquotient.
 \begin{enumerate}
  \item The Tits group $^2F_4(2)$ with $p=3$;
  \item the largest Janko group $J_4$ with $p=3$;
  \item the Rudvalis group $Ru$ with $p=3$; 
  \item the Thompson group $Th$ with $p=5$; 
  \item an exceptional simple group of Lie type denoted $^2F_4(q)$, where $q=2^{2n+1}$, $p=3$, and $9$ does not divide $q+1$.
  \end{enumerate}
\end{thm}

Recently Tiep and Navarro established that when $p\geq 7$, if all $p$-elements have class size coprime to $p$, the Sylow $p$-subgroup is abelian \cite{Tiep}.
This result is encompassed in the proof of our main theorem.
Our proof was established independantly to the work of Tiep and Navarro, and in addition to their result we consider the cases $p=3$ and $5$.

We note that the list of simple $cl_p$ groups in Theorem~\ref{MainTheorem} is the same as the list given by Tiep and Navarro in Theorem A \cite{Tiep2}; where they classified groups in which any two non-trivial $p$-elements are conjugate. Note that if any two $p$-elements are conjugate in a group and the group has non-abelian Sylow $p$-subgroups, then the group will be $cl_p$. 
However for any $cl_p$ group $G$, if $A$ is an abelian group with order divisible by $p$, then $G\times A$ is also $cl_p$ but all $p$-elements will not be conjugate; thus a classification of $cl_p$ groups can not be the same as the classification given by Theorem A \cite{Tiep2}.

Our proof of the main theorem differs from the proof by Tiep and Navarro as Tiep and Navarro established that a minimal $cl_p$ group must be a quasi-simple group, and then they only study such groups for $p\geq 7$; while in our proof, we show that a minimal $cl_p$ group must in fact be simple, and then study which simple groups can be $cl_p$. 
The advantage in our approach is that it allows us to deal with the cases  $p=3$ and $5$ only when a minimal example is a perfect central extension of a simple group which has an abelian Sylow $p$-subgroup, instead of having to study all the quasi-simple groups as perfect central extensions of a simple group with no information given about the Sylow $p$-subgroup.

Our proof of the main theorem makes use of the Classification of Finite Simple Groups.
In addition we make use of the observation that if $p$ does not divide the class size of $x$ in $G$, this is equivalent to saying that the centraliser of $x$ in $G$, denoted $C_G(x)$, contains a Sylow $p$-subgroup of $G$.

Due to Camina and Herzog's result we shall assume that $p$ is an odd prime for the rest of this paper.

At times we will make use of the notation $|X|_p$ to mean the $p$-part of the size of the set $X$. The rest of the notation used throughout is standard.

\section{Preliminary results}
We begin proving some small results which will be used during this paper. 

Proposition~\ref{exp} will be used when we study the classical groups; this is due to the fact that most of the Sylow $p$-subgroups can be realised as wreath products, and this result will allow us to exclude these groups.
Proposition~\ref{ClassLift} and Lemma~\ref{ClassSpilt} are mainly used to deal with the prime $3$, as in these cases we will need to study the Schur cover of the simple group.
Propositions~\ref{Auto2} and ~\ref{SylFromNorm} will be used in conjunction in Section 4 to find the structure of a minimal $cl_p$ group.

\begin{lm}\label{cent}
Let $G=H\wr A$ with $A$ an abelian transitive permutation group on the set $\{1,\dots, n\}$. If $g\in Z(G)$, then there exists an element $h\in Z(H)$ such that $g=(h,h,\dots,h)\in H\wr A$
\begin{proof}
Pick $g\in Z(G)$ and assume $g= (h_1, \dots, h_n)a_1$ for $h_j\in H$ and $a_1\in A$. 
As $g\in Z(G)$ then $ga=ag$ for all $a\in A$. 
In particular we have that $ga= (h_1, \dots, h_n)a_1a$ and $ag= (h_{a(1)},\dots,h_{a(n)})aa_1$. As $A$ was abelian $a$ and $a_1$ commute; hence for any permutation $a$ we have $(h_1,\dots,h_n)=(h_{a(1)},\dots,h_{a(n)})$.
By the transitivity of $A$ it follows that $h_1=h_j$ for all $j$.

Choose an element $k=(k_1,\dots, k_1)$, for $k_1$ in $H$. As $gk=kg$ this implies that $h_1k_1=k_1h_1$ for all $k_1\in H$; thus $h_1$ must in fact lie in $Z(H)$.

If $a_1$ is non trivial, then say $a_1(i)=j\not= i$. Pick $k=(k_1,\dots, k_n)$ such that $k_i\not=k_j$ and $k_r=1_H$ otherwise. Setting $gk=kg$ implies that $k_ih_1=h_1k_j$ and so $k_i=k_j$, thus $a_1$ must be trivial.
\end{proof}
\end{lm}

\begin{prop}\label{exp}
Suppose $G=H\wr A$, where $A$ is an abelian transitive permutation group on $n$ elements, then $\exp(Z(G))=\exp(Z(H))\leq \exp(H)$.
 
 If, in addition there exists $h\in H$ with the order of $h$ equal to $\exp(H)$ then $\exp(H)<\exp(G)$; it follows that $\exp(Z(G))<\exp(G)$.
\begin{proof}
The first inequality follows from Lemma~\ref{cent}.

For the second inequality, let $h\in H$ be an element of order $s=\exp(H)$. Take $1_H$ to denote the identity of $H$, and pick a non-trivial permutation $a$ which does not fix 1. 
For $g=(h,1_H,\dots,1_H)a$ it follows that $g^{s} =(h^{r_1},\dots, h^{r_n})a^s\not= 1$ with each $r_i< s$; it is now clear that $\exp(H)< \exp(G)$.
\end{proof}
\end{prop}

\begin{prop}\label{ClassLift}
 Let $G$ be a group, $N\lhd G$ and set $\overline{G}=G/N$. If $\overline{g}=gN$ is an element of $\overline{G}$, and $A_g:=\{h\in G|[g,h]\in N\}$ then $C_{\overline{G}}(\overline{g})=A_g/N$.
 
 In addition 
 \[
  |A_g|=|\{[g,G]\cap N\}|\cdot |C_G(g)|,
 \]
and so
 \[
  |C_{\overline{G}}(\overline{g})|=|\{[g,G]\cap N\}|\cdot |C_G(g)|/|N|.
 \]
\begin{proof}
Let $A_g$ be as above. If $[g,h]$ and $[g,k]$ are in $N$, then $[g,hk]=[g,k][g,h]^k$ and $[g,h^{-1}]=([g,h]^{-1})^{h^{-1}}$ are both in $N$. Thus $A_g\leq G$ and clearly $N\lhd A_g$.
From the definition of $A_g$ it is clear that $A_g/N=\{hN | [g,h]\in N \}=C_{\overline{G}}(\overline{g})$.

Let $C=C_G(g)$, then $[g,f_1]=[g,f_2]$ if and only if $f_2\in Cf_1$. Hence for each commutator $[g,h]$ which lies in $N$, the subgroup $A_g$ contains the set of elements $Ch$; so it follows that $|A_g|=|C|\cdot |\{[g,G]\cap N\}|$. 

\end{proof}
\end{prop}

\begin{lm}\label{Spilt}\label{ClassSpilt}
 Let $G$ be a finite group and $\overline{G}\cong G/Z(G)$ with the canonical homomorphism $\pi:G \rightarrow \overline{G}$. For any $x\in G$, then $\pi^{-1}(\overline{x}^{\overline{G}})$ is a union of $k$ distinct conjugacy classes, for some integer $k$.  

 It follows that $|C_G(x)|=|C_{\overline{G}}(\overline{x})|k$.
 \begin{proof}
 \underline{Claim: } $\pi^{-1}(\overline{x}^{\overline{G}})=\bigcup_{z\in Z(G)}(xz)^G=\bigsqcup_{i=1}^k(xz_i)^G$ for some $z_1,\dots,z_k\in Z$.\newline

Let $\overline{y}=yZ(G)$, then $\overline{y}=\overline{x}^{\overline{g}}$ if and only $y=x^gz=(xz)^g$ for some $z\in Z(G)$; which gives the first equality.
The second equality follows by choosing $z_1,\dots, z_k\in Z(G)$ so that each class is distinct.

For the final statement, our claim implies that $|\pi^{-1}(\overline{x}^{\overline{G}})|=k|x^G|$; but also $|\pi^{-1}(\overline{x}^{\overline{G}})|=|Z(G)||\overline{x}^{\overline{G}}|$, and so $|x^G|=|\overline{x}^{\overline{G}}|\frac{|Z(G)|}{k}$. The result now follows by applying the Orbit-Stabilizer theorem.
 
 \end{proof}
\end{lm}

\begin{prop}\label{Auto2}
 Let $G$ be a group and $N_1,\dots, N_t$ distinct normal subgroups of $G$, such that $[N_i,N_j]=1$ for all $i\not= j$. Take $\phi$ the canonical homomorphism $\phi:G\rightarrow $Aut$(N_1)\times\dots\times $Aut$(N_t)$ defined componentwise by $\phi_i:G\rightarrow$ Aut$(N_i)$ under the conjugation action. Let $\Phi$ be the composite map defined by\newline
\[
  \xymatrix{
   G \ar[rrr]^\phi \ar[rrrd]_\Phi & & 
   &\text{Aut}(N_1)\times\dots\times \text{Aut}(N_t) \ar@{->>}[d]\\
   & & & \text{Out}(N_1)\times\dots\times \text{Out}(N_t)}.
\]
Then, Ker$(\Phi)=$Ker$(\phi)N_1N_2\dots N_k$ and Ker$(\phi)=C_G(N_1)\cap\dots\cap C_G(N_k)$.
\begin{proof}
 As $[N_i,N_j]=1$ for all $i\not= j$ then the $N_i$ commute.
 
 Pick $g\in$ Ker$(\phi) N_1N_2\dots N_t$, then $g=kh_1\dots h_t$ for $h_i\in {N_i}$ and $k\in$ Ker$(\phi)$. For all $x_i\in N_i$ it follows that $x_i^g=x_i^{h_i}$; so $\phi_i(g)\in$ Inn$(N_i)$ for all $i$ and $\Phi(g)=1$.
 
 Conversely, let $g\in$ Ker$(\Phi)$; for all $i$ there exists $h_i\in N_i$ such that $x_i^g=x_i^{h_i}$ for all $x_i\in N_i$. Take $g'=h_1\dots h_t$, then for all $x_i\in N_i$, $x_i^{g(g')^{-1}}=x_i$ and hence $g(g')^{-1}\in$ Ker$(\phi)$; so $g\in$ Ker$(\phi)N_1\dots N_t$.
 
 Finally it remains to find Ker$(\phi)$. Now $\phi(g)=1$ implies $x_i^g=x_i$ for all $x_i\in N_i$ and so $g\in C_G(N_1)\cap \dots \cap C_G(N_t)$. 
 The converse is clear.
\end{proof}
\end{prop}

We also want to recall the following result. 
\begin{lm}\label{SylSubs}\label{NormalSylSubs}\cite[Proposition 6.4.7]{Huppert}
 Let $G$ be a finite group and $H,K\leq G$ such that $G=HK$ and $p$ a prime.
 In particular if $N$ and $K$ are normal subgroups of $G$ such that $G=NK$, and $P$ is a Sylow $p$-subgroup of $G$; then $P=(P\cap N)(P\cap K)$.
\end{lm}

\begin{prop}\label{SylFromNorm}
 Let $G$ be a group with 3 distinct normal subgroups $N_1,N_2,N_3$ such that $N_1N_2N_3\not=N_iN_j$ for $i\not=j$ in $\{1,2,3\}$. If all proper subgroups of $N_1N_2N_3$ have abelian Sylow $p$-subgroups, then the Sylow $p$-subgroups of $N_1N_2N_3$ are abelian.
 \begin{proof}
  Let $P$ be a Sylow $p$-subgroup of $N_1N_2N_3$; by Lemma~\ref{NormalSylSubs}, $P=(P\cap N_1N_2)(P\cap N_3)$. Each $P\cap N_iN_j$ for $i\not=j$ and $P\cap N_i$ is abelian by assumption. As $P\cap N_3$ and $P\cap N_2\leq P\cap N_2N_3$, then $P\cap N_3$ commutes with $P\cap N_2$ and similarly they commute with $P\cap N_1$.
  
  It remains to show for $g\in P\cap N_1N_2$ and $h\in P\cap N_3$, then $gh=hg$. 
  As $P\cap N_1N_2$ is a Sylow $p$-subgroup of $N_1N_2$, it follows $P\cap N_1N_2=(P\cap N_1)(P\cap N_2)$.
  Hence $g=h_1h_2$ such that $h_i\in P\cap N_i$ and $hh_i=h_ih$. It follows that $gh=hg$ and $P$ is abelian.
 \end{proof}
\end{prop}

\section{Simple groups}
We begin this section by noting the following result.
\begin{lm}
 Let $G$ be a finite group and $p$ a prime; if the class size of a $p$-element $x$ is not divisible by $p$ then there exists a Sylow $p$-subgroup $P$ such that $x\in Z(P)$.
 \begin{proof}
  Choose $P$ a Sylow $p$-subgroup of $G$, such that $P\leq C_G(x)$. 
  It is an easy observation that $P=P\times \langle x\rangle$ and so $x$ must lie in $Z(P)$.
 \end{proof}
\end{lm}

Hence if we select a property of a subgroup which is preserved by conjugation, and show that $Z(P)$ and $P$ exhibit differing properties; then the group can not be $cl_p$. 
For example, if for $P$ a Sylow $p$-subgroup of $G$ then if $Z(P)$ and $P$ have differing exponents the group can not be $cl_p$; we will use this example to cover a large number of cases for the classical simple groups.

\begin{lm}\label{IndexQuot}
Let $G$ be a group, $p$ a prime and $N$ a normal subgroup of $G$. If all $p$-elements have class size not divisible by $p$, then the same holds in both $N$ and $G/N$.
\begin{proof}
This result follows by noting that the Sylow $p$-subgroups of $N$ and $G/N$ take the form $P\cap N$ and $PN/N$ respectively, for $P$ a Sylow $p$-subgroup of $G$.
\end{proof}
\end{lm}
However, given conditions upon a normal subgroup $N$ we can also lift $cl_p$ up from a normal subgroup or quotient.
\begin{prop}\label{RestrictCLp}
 Let $G$ be a group, $N$ a normal subgroup and $p$ a prime.
 \begin{enumerate}
  \item If $p$ does not divide $|G/N|$ then $N$ is $cl_p$ if and only if $G$ is $cl_p$.
  \item If $p$ does not divide $|N|$ then $G/N$ is $cl_p$ if and only if $G$ is $cl_p$.
 \end{enumerate}
\begin{proof}
For each statement there is only one implication to show by Lemma~\ref{IndexQuot}.
 \begin{enumerate}
  \item Let $g$ be a $p$-element of $G$, then $g$ lies in $N$. Choose $P$ a Sylow $p$-subgroup of $N$ such that $P\subset C_N(g)$. Then $P$ is a Sylow $p$-subgroup of $G$ and so $P\subset C_G(g)$.
  \item Let $g$ be a $p$-element of $G$, then $gN$ must lie in the centre of a Sylow $p$-subgroup $PN/N$ of $G/N$, for $P$ a Sylow $p$-subgroup of $G$. As $Z(PN/N)=Z(P)N/N$, then $g$ must lie in $Z(P)$. 
 \end{enumerate}
\end{proof}
\end{prop}

The aim of this section is to prove that the list in Theorem~\ref{MainTheorem} is the list of all possible simple $cl_p$ groups.

\begin{thm}\label{SimpleClp}
Let $G$ be a finite simple group, and let $p$ be a prime. If $G$ is a $cl_p$ group, then $G$ must be one of the following groups.
 \begin{enumerate}
  \item The Tits group $^2F_4(2)'$ with $p=3$;
  \item the largest Janko group $J_4$ with $p=3$;
  \item the Rudvalis group $Ru$ with $p=3$; 
  \item the Thompson group $Th$ with $p=5$; or
  \item an exceptional simple group of Lie type denoted $^2F_4(q)$ for $q=2^{2n+1}$, $p=3$ and where $9$ does not divide $q+1$.
 \end{enumerate}
\end{thm}

\begin{lm}
No alternating group is $cl_p$ for any prime.
 \begin{proof}
  Let $G$ be the alternating group $A_n$ for some $n$. 
  
  If $n<p^2$, then the Sylow $p$-subgroup is elementary abelian; so assume that $n\geq p^2$. 
  By choosing $g$ to be a cycle of length $p^2$ and as we have assumed $p$ is odd we have $|g^{S_n}|_p=|g^{A_n}|_p$; it is therefore enough to consider the conjugacy class in $S_n$.
  For such a $g$ we get $C_{S_n}(g)\cong C_{S_{p^2}}(g)\times S_{n-p^2}\cong C_{p^2}\times S_{n-p^2}$; by comparing group orders it follows that $p$ divides the class size of $g$. 
 \end{proof}
\end{lm}

\begin{lm}
 The only sporadic simple groups which are $cl_p$ are given by the following list.
 \begin{enumerate}
  \item The largest Janko group $J_4$ with $p=3$;
  \item the Rudvalis group $Ru$ with $p=3$; or
  \item the Thompson group $Th$ with $p=5$.
 \end{enumerate}
 \begin{proof}
  These groups are checked directly using the Atlas \cite{Atlas}.
 \end{proof}
\end{lm}

\begin{lm}\label{unipotent}
 Let $G$ be a finite simple group of Lie type over a field of characteristic $p$, then $G$ is not $cl_p$.
 \begin{proof}
  Take $G$ a simple algebraic group of adjoint type and $F:G\mapsto G$ a Frobenius map, then the corresponding finite simple group can be viewed as $S:=[G^F,G^F]$ for $G^F$ the fixed point set of the Frobenius map $F$ \cite[Table 22.1]{TestMal}.
  
  Following the notation of Luo in \cite{Betty}, let $u$ be the regular unipotent element given by $u:=\prod_{\alpha\text{ a simple root}}x_\alpha(1)$, which gives a $p$ element in the finite simple group of Lie type. 
  Then $C_S(u)\leq C_G(u)$ for $G$ the simple algebraic group; however Lou established that $C_G(u)$ is abelian \cite[Page 1145]{Betty}. Hence if the regular unipotent element contains a Sylow $p$-subgroup in its centraliser then in fact the Sylow $p$-subgroup must be abelian. 
 \end{proof}
\end{lm}
\subsection{The classical simple groups of Lie type}
To study the classical simple groups over the finite field $\F_q$ when $p$ does not equal the field characteristic, we can make use of Proposition~\ref{RestrictCLp}. 
In the symplectic and orthogonal groups, the restriction from the classical group to the simple group requires us to divide the order by a power of 2 at each stage \cite[Sections 3.7 and 3.8]{Fsg}; therefore it is enough to show the classical group is not $cl_p$. 
For the remaining two familes, the projective special linear and the projective special unitary groups; we can also prove the case for $p$ not dividing the matrix size $n$ and $q-1$ or $q+1$ respectively by establishing the result for the general linear group and the unitary group and applying Proposition~\ref{RestrictCLp}.

To achieve this we recall that for classical groups the Sylow $p$-subgroups in cross characteristic occur as wreath products.
\begin{thm}\cite[Page 6]{WeirClass}\label{ClassSylow}
The Sylow $p$-subgroups of the classical groups over $\F_q$ with $q$ prime to $p$ and $p$ odd are expressible as direct products of groups defined inductively by $G_i:=G_{i-1}\wr C_p$, where $G_0$ is an abelian $p$-subgroup of the classical group.
\end{thm}
We note that $P$ a Sylow $p$-subgroup of a classical group is abelian if and only if $P\cong G_0\times\dots\times G_0$. 

\begin{cor}\label{ClassicalSylowClp}
 Let $G$ be a finite group and $p$ a prime. If a Sylow $p$-subgroup of $G$ is isomorphic to a Sylow $p$-subgroup of a classical group then $G$ is not $cl_p$. 
\begin{proof}
 By applying Proposition~\ref{exp} to Theorem~\ref{ClassSylow} it becomes clear that if a classical group has a non-abelian Sylow $p$-subgroup $P$, then $exp(P)>exp(Z(P))$; and so in fact can not be a $cl_p$ group.
\end{proof}
\end{cor}

It only remains to show that the classical simple groups $PSL_n(q)$ and $PSU_n(q)$ are not $cl_p$ when the Sylow $p$-subgroup of the simple group is not isomorphic to a Sylow $p$-subgroup of the classical group. In particular we will show that unless $p$ divides the matrix size, then we can apply Corollary~\ref{ClassicalSylowClp}.

\paragraph{\underline{The projective linear group}}

We split the problem into two cases, $p$ divides $q-1$ but not $n$, and then $p$ divides $n$ as well.

If $p$ does not divide $n$, by Proposition~\ref{RestrictCLp} it is enough to show that $SL_n(q)$ is not $cl_p$. 
For $SL_n(q)$ we can obtain the Sylow $p$-subgroup from $GL_n(q)$ using the technique outlined in Weir's paper \cite[Section 2]{WeirClass}. 
We write $n=a_0+a_1p+\dots+a_sp^s$ where each $a_i<p$; then a Sylow $p$-subgroup of $GL_n(q)$ is given by $P\cong G_0^{a_0}\times \dots\times G_s^{a_s}$, for $G_i$ as in Theorem~\ref{ClassSylow}. 
As $p$ does not divide $n$ then $a_0\not=0$, and it can be shown that $P\cap SL_n(q)$ is isomorphic to $P'=G_0^{a_0-1}\times \dots\times G_s^{a_s}$, as $G_0$ is isomorphic to the Sylow $p$-subgroup of $\F_q^*$.
However, $P'$ is a Sylow $p$-subgroup of $GL_{n-1}(q)$; thus as $p$ divides $q-1$ but not $n$ it follows a Sylow $p$-subgroup of $PSL_n(q)$ is isomorphic to a Sylow $p$-subgroup of $GL_{n-1}(q)$ and so $PSL_n(q)$ is not $cl_p$ by Corollary~\ref{ClassicalSylowClp}.

Finally assume that $p|gcd(n,q-1)$, so we can assume $n\geq 3$.
Let $g=$ diag$(\zeta,\zeta^{-1},1,\dots,1)$ for $\zeta$ a $p$-element in $\F_q$ and choose $h=(h_{ij})_{1\leq i,j\leq n}$ an element in $C_{PSL_n(q)}(gZ(SL_n(q)))$.
If $n\geq 4$ it can be shown that $C_{PSL_n(q)}(gZ(SL_n(q)))=C_{SL_n(q)}(g)/Z(SL_n(q))$ and in addition $C_{SL_n(q)}(g)\cong GL_{n-2}(q)\times C_{q-1}$, where $C_{q-1}$ is considered to be an isomorphic copy of the field $\F_q^*$.
The index of $C_{PSL_n(q)}(gZ(SL_n(q)))$ in $PSL_n(q)$ is $(q^{n-1}+\dots +1)(q^{n-2}+\dots +1)$.  
The following lemma shows that the element $gZ(S_n(q))$ has class size divisible by $p$.
\begin{lm}
 Let $p$ be a prime and $q$ a prime power such that $p|q-1$ then $p|q^{k-1}+\dots +1$ if and only if $p|k$. 
\end{lm}
The proof of this lemma is a straight forward inductive argument using that when a prime $p$ divides $ab$ then $p$ divides $a$ or $b$.

If $n=3$ then we must have $p=3$.
If $9$ does not divide $q-1$, then it is easy to see that the Sylow $3$-subgroup of $PSL_3(q)$ has order 9 and so is abelian. 
Thus we can assume that the field $\F_q$ has an element $\zeta$ of order $9$, and let $g=$ diag$(\zeta,\zeta^{-1},1)$.
As $\zeta$ has order $9$, it can be shown that $C_{PSL_3(q)}(gZ(SL_3(q))=C_{SL_3(q)}(g)/Z(SL_3(q))$ and in addition $C_{SL_3(q)}(g)\cong C_{q-1}^2$.
The index of $C_{PSL_3(q)}(gZ(SL_3(q)))$ in $PSL_3(q)$ is $(q+1)(q^2+q+1)$; by the previous lemma we see that $3$ divides the class size of this element.

\paragraph{\underline{The projective unitary group}}
The argument for $PSU_n(q)$ will be similar to the case of $PSL_n(q)$.

The problem splits into two cases, $p$ divides $q+1$ but not $n$, and then $p$ divides $n$ as well.

If $p$ does not divide $n$, by Proposition~\ref{RestrictCLp} it is enough to show that $SU_n(q)$ is not $cl_p$. 
As $p|q^2-1$ it follows by comparing group orders that a Sylow $p$-subgroup of $SU_n(q)$ is a Sylow $p$-subgroup of $SL_n(q^2)$. 
As $p$ does not divide $n$, it follows that a Sylow $p$-subgroup of $SL_n(q^2)$ is isomorphic to a Sylow $p$-subgroup of $GL_{n-1}(q^2)$ by the same argument as in the projective linear case.
Thus $PSU_n(q)$ is not a $cl_p$ group by Corollary~\ref{ClassicalSylowClp}.

Finally assume that $p|gcd(n,q+1)$, so we can assume $n\geq 3$.
Let $g=$ diag$(\zeta,\zeta^{-1},1,\dots,1)$ for $\zeta$ a $p$-element in $\F_{q^2}$.
If $n\geq 4$ by the same argument as in the linear case, $C_{PSU_n(q)}(gZ(SU_n(q)))=C_{SU_n(q)}(g)/Z(SU_n(q))$ and $C_{SU_n(q)}(g)\cong GU_{n-2}(q)\times C_{q+1}$, where $C_{q+1}$ is considered to be the cyclic subgroup of order $q+1$ in the field $\F_{q^2}^*$.
The index of $C_{PSU_n(q)}(gZ(SU_n(q)))$ in $PSU_n(q)$ is $(q^{n-1}-q^{n-2}+\dots-q+1)(q^{n-3}+q^{n+5}+\dots +q^2+1)(q-1)$ if $n$ is odd and $(q^{n-2}+q^{n-4}+\dots +q+1)(q^{n-2}-q^{n-3}+\dots -q+1)(q-1)$ if $n$ is even.
The following lemma shows that the element $gZ(SU_n(q))$ has class size divisible by $p$.

\begin{lm}
 Let $p$ be a prime and $q$ a prime power such that $p|q+1$ then:
\begin{enumerate}
 \item for $k$ odd, $p$ divides $q^{k-1}-q^{k-2}+\dots -q+1$ if and only if $p$ divides $k$,
 \item for $k$ even, $p$ divides $q^{k-2}+q^{k-4}+\dots+q^2+1$ if and only $p$ divides $\frac{k}{2}$.
\end{enumerate}
\end{lm}
The proof of this result follows the same idea as the proof of Lemma 3.10.

If $n=3$ then we must have $p=3$.
If $9$ does not divide $q+1$, then it is easy to see that the Sylow $3$-subgroup of $PSU_3(q)$ has order 9 and so is abelian. 
Thus we can assume that the field $\F_{q^2}$ has an element $\zeta$ of order $9$, and let $g=$ diag$(\zeta,\zeta^{-1},1)$.
As $\zeta$ has order $9$, it can be shown that $C_{PSU_3(q)}(gZ(SU_3(q))=C_{SU_3(q)}(g)/Z(SU_3(q))$ and in addition $C_{SU_3(q)}(g)\cong C_{q+1}^2$.
The index of $C_{PSU_3(q)}(gZ(SU_3(q)))$ in $PSU_3(q)$ is $(q-1)(q^2-q+1)$; by the previous lemma we see that $3$ divides the class size of this element.

\begin{prop}
 Let $G$ be a classical simple group and $p$ a prime. Then $G$ is not a $cl_p$ group.
\end{prop}

\subsection{The exceptional simple groups of Lie type}
In this section we prove that if an exceptional simple group of Lie type is $cl_p$, then $p=3$ and the only such simple groups are $^2F_4(2)'$ and $^2F_4(q)$ for $q=2^{2n+1}$ such that $9$ does not divide $q+1$.

First we exclude the group $^2F_4(2)'$. By using the Atlas \cite{Atlas} to look at class sizes and subgroup structure, we can see this group is $cl_p$ only for $p=3$.

To study the rest of these simple groups we use \cite[Section E.II, Corollary 5.19]{SAG} which tells us that if $p$ is not the underlying field characterisctic and if $p$ does not divide the order of the corresponding Weyl group; then the Sylow $p$-subgroups of the exceptional group of Lie type are abelian.
The only cases left to consider are given below.\newline
\begin{center}

\begin{tabular}{c|c|c|c|c|c}
Group & $^2B_2(2^{2n+1})$ & $G_2(q)$ & $^2G_2(3^{2n+1})$ & $^3D_4(q)$ & $F_4(q)$\\
\hline
Primes & $2$ & $2,3$ & $2,3$ & $2,3$ & $2,3$\\
\end{tabular}

\begin{tabular}{c|c|c|c|c|c}
 Group & $^2F_4(2^{2n+1})$ & $E_6(q)$ & $^2E_6(q)$ & $E_7(q)$ & $E_8(q)$\\
 \hline
 Primes  & $2,3$ & $2,3,5$ & $2,3,5$ & $2,3,5,7$ & $2,3,5,7$\\ 
\end{tabular}
\end{center}

As we assume $p$ is odd and $p$ is not the field characteristic, then we do not need to consider either $^2B_2(2^{2n+1})$ or $^2G_2(3^{2n+1})$.

For the other cases we assume $p$ is not the field characteristic by Lemma~\ref{unipotent}.
We either use tables to find an element in the group with class size divisible by $p$; or we use tables about maximal subgroups and reduce the argument down to the classical case and apply Corollary~\ref{ClassicalSylowClp}.

We can express the order of an exceptional simple group of Lie type $G$ as $q^a\cdot \prod_i \Phi^{b_i}_{m_i}(q)$ for $q$ a power of $p$ , and where $a,b_i$ and $m_i$ are positive integers, with the $m_i$ being pairwise distinct and $\Phi_m(q)$ denotes the $m^{th}$ cyclotomic polynomial \cite[4.10.1]{CFSG}; from now on we shall write $\Phi_m$ instead of $\Phi_m(q)$. 
In addition if $p$ only divides one of the cyclotomic polynomials, then the Sylow $p$-subgroup is abelian \cite[4.10.2]{GorLyon}.

As the order of an exceptional simple group is divisible only by $\Phi_m$ if $1\leq m\leq 30$, for this section we shall only refer to cyclotomic polynomials where $1\leq m\leq 30$.

For the prime $p=7$ we only have the groups $E_7(q)$ and $E_8(q)$; these cases have been considered by Navarro and Tiep in \cite{Tiep}, however we include the proof for completeness.
\begin{lm}
 No exceptional simple group of Lie type is $cl_7$.
\begin{proof}
If $q\not\equiv 1,-1$ modulo $7$, then $7$ only divides one of the cyclotomic polynomials and so the Sylow $7$-subgroup of $E_7$ and $E_8$ will be abelian. 
Thus assume $q\equiv 1$ or $-1$ modulo 7.
If $q\equiv 1$ modulo 7 then 7 divides $\Phi_1$ and $\Phi_7$; while if $q\equiv -1$ modulo 7 then 7 divides $\Phi_2$ and $\Phi_{14}$.

\underline{\textbf{The group $E_7(q)$}}.
Suppose $q\equiv 1$ modulo 7, in $E_7(q)$ we find a maximal subgroup of order divisible by $\Phi_1^7 \Phi_7$; by \cite[Table 5.1]{LieSaxSeit} the subgroup isomorphic to $PSL_8(q)$ has order divisible by $\Phi_1^7\Phi_7$.
Additionally $7$ does not divide $8$ and so a Sylow 7-subgroup of $PSL_8(q)$ is isomorphic to a Sylow 7-subgroup of $GL_7(q)$; thus $E_7(q)$ is not $cl_7$ by Corollary~\ref{ClassicalSylowClp}.

Similary for $q\equiv -1$ modulo 7, then the order of $E_7(q)$ has the factor $\Phi_2^7 \Phi_{14}$; by \cite[Table 5.1]{LieSaxSeit} we consider the subgroup isomorphic $PSU_8(q)$.
In same way as the previous case, a Sylow 7-subgroup of $PSU_8(q)$ is isomorphic to a Sylow 7-subgroup of $GL_7(q^2)$ and by Corollary~\ref{ClassicalSylowClp} the group $E_7(q)$ is not $cl_7$.

\underline{\textbf{The group $E_8(q)$}}.
As $E_8(q)$ has the subgroup $PSL_2(q)\times E_7(q)$ \cite[Table 5.1]{LieSaxSeit}, we can use the case for $E_7(q)$ to show that this group is not $cl_7$. 
The order of this group is divisible by both $\Phi_1^8 \Phi_7$ and $\Phi_2^8 \Phi_{14}$ and using the case for $E_7(q)$ we get that the Sylow $7$-subgroup of $E_8(q)$ must be the direct product of Sylow $p$-subgroups of two classical groups. Thus by Corollary~\ref{ClassicalSylowClp} the group $E_8(q)$ is not $cl_7$.
\end{proof}
\end{lm}

For $p=5$ the ideas are similar as for $p=7$.
However we can not exclude any cases for $q$ modulo 5. In fact if $q\equiv 1$ modulo $5$ then $5$ divides $\Phi_1$ and $\Phi_5$, if $q\equiv -1$ then $5$ divides $\Phi_2$ and $\Phi_{10}$ and if $q\equiv \pm 2$ modulo $5$ then $5$ divides $\Phi_4$ and $\Phi_{20}$.
\begin{lm}
No exceptional simple group of Lie type is $cl_5$.
 \begin{proof}
 For $p=5$, we have to consider the groups $E_6(q)$, $^2E_6(q)$, $E_7(q)$ and $E_8(q)$. The table below gives a list of maximal subgroups for the groups $E_6(q)$, $^2E_6(q)$ and $E_7(q)$, in the cases where 5 divides two of the cyclotomic polynomials.
 
 \begin{tabular}{c|c|c|c}
 Group & Cyclotomic factors of $|G|$ divisible by 5 & $q$ modulo 5 & Maximal subgroup\\
 \hline
  & & & \\
 $E_6(q)$ & $\Phi_1^6\Phi_2^4\Phi_4^2\Phi_5$ & $q\equiv 1$ mod 5 & $PSL_2(q)\times PSL_6(q)$\\
\hline
 $^2E_6(q)$ & $\Phi_1^4\Phi_2^6\Phi_4^2\Phi_{10}$ & $q\equiv -1$ mod 5 & $PSL_2(q)\times PSU_6(q)$\\
\hline
 $E_7(q)$ & $\Phi_1^7\Phi_2^7\Phi_4^2\Phi_5\Phi_{10}$ & $q\equiv 1$ mod 5 & $PSL_8(q)$\\
  &  & $q\equiv -1$ mod 5 & $PSU_8(q)$\\
 \end{tabular}
 
 Using the table above, the Sylow $5$-subgroup is isomorphic to a Sylow 5-subgroup of a classical group and so by the same arguments as we used for $p=7$ and by applying Corollary~\ref{ClassicalSylowClp} we obtain that these groups are not $cl_5$.
 
 This leaves us with $E_8(q)$, in this case we have the factor $\Phi_1^8\Phi_2^8\Phi_4^4\Phi_5^2\Phi_{10}^2\Phi_{20}$; however we can no longer use a maximal subgroup argument as we find that for the Sylow 5-subgroup, the exponent of the group and its centre are equal. 
 Instead consider \cite[Table 4.7.3B]{CFSG}, where a list of centralisers of elements of order 5 in $E_8(q)$ has been given. In the following table we have chosen a 5 element such that 5 divides its class size.
 
 \begin{tabular}{c|c|c}
 $q$ mod 5 & Centraliser type & The order of the 5 part of the centraliser\\
 \hline
  & & \\
 $q\equiv \epsilon=\pm1$ & $E_7(q)+ (q-\epsilon)$ & $\Phi_{\epsilon}\Phi_1^7\Phi_2^7\Phi_4^2\Phi_5\Phi_{10}$\\
 $q^2\equiv -1$ & $D_6^-(q)+(q^2+1) $ & $\Phi_1^6\Phi_2^6\Phi_4^4\Phi_5^1\Phi_{10}$\\
 \end{tabular}

\end{proof}
\end{lm}

For the prime $p=3$ we will make use of tables by L\"{u}beck \cite{Lubeck}; these tables provide a list of all elements of order 3 in exceptional simple groups of Lie type with their centraliser order. 
We have excluded $E_6(q),$ $^2E_6(q)$ and $E_7(q)$ as these groups have adjoint and simply connected versions which are distinct from the simple group \cite[Table 22.1]{TestMal} and so will need to be treated differently. 

We make the following observations, if $q\equiv 1$ modulo 3, then 3 divides $\Phi_1,\Phi_3$ and $\Phi_9$; while if $q\equiv -1$ modulo $3$ then $3$ divides $\Phi_2,\Phi_6$ and $\Phi_{18}$.

For each group in the table below, we have picked a 3-element from L\"{u}beck's table \cite{Lubeck} such that $3$ will divide the class size of the 3-element.

\begin{center}
\begin{tabular}{c|c|c}
Group & Centraliser label & Centraliser order\\
\hline
$G_2(q)$ and $q\equiv 1$ mod 3 & $\sim A_1(q) + (q-1)$ & $q\Phi_1^2\Phi_2$\\
$G_2(q)$ and $q\equiv -1$ mod 3 & $\sim A_1(q) + (q+1)$ & $q\Phi_1\Phi_2^2$\\
$^3D_4(q)$ and $q\equiv 1$ mod 3 & $A_1(q^3) + (q-1)$ & $q^3\Phi_1^2\Phi_2\Phi_3\Phi_6$\\
$^3D_4(q)$ and $q\equiv -1$ mod 3 & $A_1(q^3) + (q+1)$ & $q^3\Phi_1\Phi_2^2\Phi_3\Phi_6$\\
$F_4(q)$ and $q\equiv 1$ mod 3 & $C_3(q) + (q-1)$ & $q^9\Phi_1^4\Phi_2^3\Phi_3\Phi_4\Phi_6$\\
$F_4(q)$ and $q\equiv -1$ mod 3 & $C_3(q) + (q+1)$ & $q^9\Phi_1^3\Phi_2^4\Phi_3\Phi_4\Phi_6$\\
$E_8(q)$ and $q\equiv 1$ mod 3 & $A_8(q)$ & $q^{36}\Phi_1^8\Phi_2^4\Phi_3^3\Phi_4^2\Phi_5\Phi_6\Phi_7\Phi_8\Phi_9$\\
$E_8(q)$ and $q\equiv -1$ mod 3 & $^2A_8(q)$ & $q^{36}\Phi_1^4\Phi_2^8\Phi_3\Phi_4^2\Phi_6^3\Phi_8\Phi_{10}\Phi_{14}\Phi_{18}$\\
\end{tabular}
\end{center}

For the groups $E_6(q)$, $^2E_6(q)$ and $E_7(q)$ L\"{u}beck's tables \cite{Lubeck} give $3$-elements for $E_6(q)_{\text{ad}},E_6(q)_{\text{sc}},$ $^2E_6(q)_{\text{ad}},$ $^2E_6(q)_{\text{sc}},E_7(q)_{\text{ad}}$ and $E_7(q)_{\text{sc}}$; which are not in fact the finite simple groups of Lie type.
Instead we will use the fact that the simply connected versions $E_6(q)_{\text{sc}},$ $^2E_6(q)_{\text{sc}}$ and $E_7(q)_{\text{sc}}$ give the Schur cover of the simple group \cite[Page 209]{CFSG}. In this situation we have the following cases.

The information on the Schur cover of these simple groups is taken from \cite[Table 8.5]{Schur}; while the centraliser information has been taken from \cite{Lubeck}.

\paragraph{\underline{$q=1$ modulo 3}}

\begin{center}
\begin{tabular}{c|c|c|c}
Group & Schur cover & Centraliser label & Centraliser order\\
\hline
$E_6(q)$ & $E_6(q)=E_6(q)_{\text{sc}}/Z$ and $|Z|=3$ & $A_5(q)+(q-1)$ & $q^{15}\Phi_1^6\Phi_2^3\Phi_3^2\Phi_4\Phi_5\Phi_6$\\ 
$^2E_6(q)$ & $^2E_6(q)=$ $^2E_6(q)_{\text{sc}}$ & $^2D_4(q)+(q^2-1)$ & $q^{12}\Phi_1^4\Phi_2^4\Phi_3\Phi_4\Phi_6\Phi_8$\\ 
$E_7(q)$ & $E_7(q)=E_7(q)_{\text{sc}}/Z$ and $|Z|=$ gcd $(2,q-1)$ & $A_6(q)+(q-1)$ & $q^{21}\Phi_1^7\Phi_2^3\Phi_3^2\Phi_4\Phi_5\Phi_6\Phi_7$\\ 
\end{tabular}
\end{center}

\paragraph{\underline{$q=-1$ modulo 3}}

\begin{center}
\begin{tabular}{c|c|c|c}
Group & Schur cover & Centraliser label & Centraliser order\\
\hline
$E_6(q)$ & $E_6(q)=E_6(q)_{\text{sc}}$ & $A_5(q)+(q+1)$ & $q^{15}\Phi_1^5\Phi_2^4\Phi_3^2\Phi_4\Phi_5\Phi_6$\\ 
$^2E_6(q)$ & $^2E_6(q)=$ $^2E_6(q)_{\text{sc}}/Z$ and $|Z|=3$ & $D_4(q)+(q^2+2q+1)$ & $q^{12}\Phi_1^4\Phi_2^6\Phi_3\Phi_4^2\Phi_6$\\ 
$E_7(q)$ & $E_7(q)=E_7(q)_{\text{sc}}/Z$ and $|Z|=$ gcd $(2,q-1)$ & $^2A_6(q)+(q+1)$ & $q^{21}\Phi_1^3\Phi_2^7\Phi_3\Phi_4\Phi_6^2\Phi_{10}\Phi_{14}$\\ 
\end{tabular}
\end{center}

As $E_7(q)$ is either the simply connected group or a quotient by the centre $Z$ of order 2; then by Proposition~\ref{ClassLift}, as $3$ divides the class size of a $3$-element in $E_7(q)_{\text{sc}}$; then $3$ will divide the class size of a $3$-element in $E_7(q)_{\text{sc}}/Z=E_7(q)$.

In the cases for $E_6(q)$ and $^2E_6(q)$ such that the simple group equals the simply connected group, then we can choose a 3-element such that $3$ divides its class size.
Finally we consider the case when the simply connected group is not equal to the simple group. 
In this case the simple group is a quotient of the simply connected version by a subgroup of order $3$. 
The centraliser in the above tables was chosen so that $9$ divides the class size of the $3$-element in the simply connected group and so by Proposition~\ref{ClassLift} (for $G$ the simply connected version and $N$ the center of $G$) we get that $3$ will divide the class size of the corresponding $3$-element in the simple group. 

This leaves the final case $^2F_4(q)$ to study, where $q=2^{2n+1}$. 
This case will use both arguments that we used so far.

The group $^2F_4(2^{2n+1})$ has order $q^{12}\Phi_{12}\Phi_{6}\Phi_4^2\Phi_3\Phi_2^2\Phi_1^2$.
As $q=2^{2n+1}$, it follows that $q\equiv -1$ mod 3 and thus the 3-part of the order of this group is in $\Phi_2^2\Phi_6$.

Using the maximal subgroups table \cite[Table 5.1]{LieSaxSeit}, we find the subgroup $SU_3(q)$ which is divisible by $\Phi_2^2\Phi_6$.
(Note as 3 divides the matrix size we can not use Corollary~\ref{ClassicalSylowClp} as the Sylow 3-subgroup can no longer be seen as being isomorphic to the Sylow 3-subgroup of a classical group.) 

As 3 divides $q+1$ then 3 divides $q^2-1$ but not $q-1$; so by Weir's paper a Sylow 3-subgroup of $SU_3(q)$ will be a Sylow $3$-subgroup of $SL_3(q^2)$ \cite[Section 4]{WeirClass}.
In this case we have that a Sylow 3-subgroup of $SU_3(q)$ can be chosen to be $P=\{$ diag$(x,y,(xy)^{-1})\cdot a| x,y\in [\F_{q^2}^{*}]_3$ and $  a\in C_3\}$. 
In the same way as we proved Lemma~\ref{cent} it can be shown that $Z(P)=\{$ diag $(x,x,x)$ such that $x^3=1\}$ and we conclude that exp($Z(P))=3$.

We have two cases. If $9$ does not divide $q+1$, then exp$(P)=3$; from L\"{u}becks table \cite{Lubeck} there is a unique conjugacy class of elements of order 3, and thus the group is $cl_3$.
If $9$ divides $q+1$ then for $x$ in $\F_{q^2}$ of order 9, the element diag$(x,x^{-1},1)$ will also have order 9; thus exp$(P)>$ exp$(Z(P))$ and the group can not be $cl_3$.

\begin{lm}
 If $G$ is an exceptional simple group of Lie type which is $cl_3$, then $G$ is either $^2F_4(2)'$ or $^2F_4(q)$ for $q=2^{2n+1}$ such that $9$ does not divide $q+1$.
\end{lm}

This completes the proof of Theorem~\ref{SimpleClp}.
\section{The Main Theorem}
The aim of this section is to classify all minimal $cl_p$ groups, where we call a group $G$ a minimal $cl_p$ group if $G$ is $cl_p$  but no normal subgroup or quotient group of $G$ is $cl_p$.

\begin{thm}
 Let $G$ be a minimal $cl_p$ group, then $G$ is a simple group.
\end{thm}
This result has the following corollary.
\begin{cor}\label{SubquotientsClp}
 Any $cl_p$ group must contain a simple $cl_p$ group as a subquotient.
\end{cor}
We shall assume throughout that $G$ is a minimal $cl_p$ group which is not simple and derive a contradiction.

\paragraph{\underline{Step 1: If $N\lhd G$ such that $1\not=N\not=G$, then $N$ and $G/N$ have abelian Sylow subgroups.}}
By Lemma~\ref{IndexQuot} all $p$-elements in $N$ and $G/N$ have class size not divisible by $p$, but by the minimality of $G$ they cannot be $cl_p$ groups.
\paragraph{\underline{Step 2: $O_{p'}(G)=1$ and $O^{p'}(G)=G$.}}
By considering Proposition~\ref{RestrictCLp} it is clear that $O_{p'}(G)=1$ and $O^{p'}(G)=G$ as otherwise a quotient of $G$ or a normal subgroup would be $cl_p$.

\paragraph{\underline{Step 3: $O_p(G)=Z(G)$.}}

Let $P$ be a Sylow $p$-subgroup of $G$ and $g\in P$; then $O_p(G)\subseteq C_G(g)$. Thus $O_p(G)\subseteq Z(P)$, or in other words $P\subset C_G(O_p(G))$.
Consider $C_G(O_p(G))\lhd G$ then the quotient group $G/C_G(O_p(G))$ is a $p'$-group and therefore trivial by Step 2, thus $O_p(G)\leq Z(G)$.

Let $Q$ be a Sylow-$q$ subgroup of $Z(G)$ for a prime $q$; as $Q$ is characteristic in $Z(G)\lhd G$ then $Q\lhd G$ and $q=p$ by Step 2.

\paragraph{\underline{Step 4: $G$ is perfect.}}
The group $G/G'$ is abelian and if $p$ divides its order then there exists a normal subgroup $K/G'$ of index $p$, with $G'\leq K\lhd G$. 
Pick $x$ a $p$-element in $G\backslash K$, then $C_G(x)$ contains $P$ a Sylow $p$-subgroup. 
Now $P\cap K$ is a Sylow $p$-subgroup of $K$ and so abelian; however this would mean that $Q:=(P\cap K)\langle x\rangle$ must be an abelian Sylow $p$-subgroup of $G$.

\paragraph{}
From now on let $\overline{G}$ denote the quotient group $G/O_p(G)$.

\paragraph{\underline{Step 5: $O_{p'}(\overline{G})=1$.}}

Set $\overline{K}=O_{p'}(\overline{G})$, then $\overline{K}=K/O_p(G)$ for $K=O_{p,p'}(G)\leq G$. By the Schur-Zassenhaus Theorem \cite[Theorem 18.1]{Asch} there exists a subgroup $H$ of $K$ such that $K=O_p(G)H$ and $H\cap O_p(G)=1$. As $O_p(G)= Z(G)$ it follows that $K=O_p(G)\times H$.

As every element in $K\backslash H$ has order divisible by $p$, then $H$ is the set of elements with $p'$-order. Thus $H$ is a characteristic subgroup of $K$; and it follows that $H$ is a normal $p'$-subgroup of $G$ and so $K=O_p(G)$ by Step 2.

\paragraph{\underline{Step 6: Each minimal normal subgroup of $\overline{G}$ is non-abelian.}}

Take $\overline{N}$ a minimal normal subgroup of $\overline{G}$; then we can write $\overline{N}\cong \overline{S_1}\times \dots\times \overline{S_t}$ \cite[Section 3.3, Page 84]{Robin} such that $\overline{S_i}\cong \overline{S_1}$ for all $i$ and $\overline{S_1}$ is simple. 

If $\overline{N}$ was abelian, each $S_i$ must be abelian. The only abelian simple groups are cyclic groups of prime order, hence $\overline{N}$ would be a normal $p$-subgroup, contradicting the fact we have quotiented out by $O_p(G)$. 

\paragraph{\underline{Step 7: Each minimal normal subgroup of $\overline{G}$ is simple.}}

$G$ acts on $\overline{N}\cong \overline{S_1}\times \dots\times \overline{S_t}$ by conjugation; let $g\in G$ and consider $\overline{S_i}^g$.
Given $\overline{h}\in \overline{N}\lhd \overline{G}$, then $g\overline{h}=\overline{h'}g$ for some $\overline{h'}\in \overline{N}$. As $\overline{S_i}\lhd \overline{N}$ it follows that $(\overline{S_i}^g)^{\overline{h}}=(\overline{S_i}^{\overline{h'}})^g=\overline{S_i}^g$ hence $\overline{S_i}^g\lhd \overline{N}$. 
As each of the $\overline{S_i}$ are non-abelian and simple; then if $\overline{S_i}^g\cap\overline{S_j}=1$, then for each element $(g_1,\dots,g_t)$ in $\overline{S_i}^g$, then for $(1,\dots,h,1,\dots 1)$ in $S_j$ we have that $[(g_1,\dots,g_t),(1,\dots,h,1,\dots 1)]=1$ and so $g_j$ must lie in $Z(S_i^g)=1$. Thus there exists a $k$ such that $\overline{S_i}^g\cap \overline{S_k}\not= 1$ otherwise $\overline{S_i}^g=1$; hence the action of $G$ must permute the $\overline{S_i}$. 

Set $K:=\cap N_G(\overline{S_i})$ the kernel of this action. The aim is to show all $p$-elements lie in $K$.

Pick $x$ a $p$-element of $G$; then $C_G(x)$ contains a Sylow $p$-subgroup $P$ of $G$. 
Let $\overline{P}:=P/O_p(G)$; then as $\overline{S_i}\lhd\overline{N}$, then $\overline{P}\cap \overline{S_i}$ is a Sylow $p$-subgroup of $\overline{S_i}$.
The conjugation action of $x$ on $\overline{N}$ fixes $\overline{P}\cap\overline{S_i}$ which is a Sylow $p$-subgroup of $\overline{S_i}$ and so $x\in K$. 
Thus $G/K$ is a $p'$-group and $G=K$ by Step 2; this means $\overline{S_1}^g=\overline{S_1}$ for all $g\in G$ and so $\overline{S_1}\lhd \overline{G}$ implying $\overline{N}=\overline{S_1}$.

\paragraph{\underline{Step 8: $\overline{G}$  is a direct product of its minimal normal subgroups.}}

Let $\overline{N_1},\dots, \overline{N_t}$ be the set of all distinct minimal normal subgroups of $\overline{G}$. As each $\overline{N_i}$ is simple, then Out$(\overline{N_i})$ is soluble by Schreier's Conjecture \cite[Page 160]{Asch}. 

As the $\overline{N_i}$ are minimal, it follows that $[\overline{N_i},\overline{N_j}]=1$ for $i\not=j$; thus we can apply Proposition~\ref{Auto2} to $\overline{G}$ with $\overline{K}=$ Ker$(\Phi)$ to see that $\overline{G}/\overline{K}$ is isomorphic to a subgroup of Out$(\overline{N_1})\times\dots\times$ Out$(\overline{N_t})$.
As $G$ is perfect the quotient $\overline{G}/\overline{K}\cong G/K$ is perfect and soluble so must be trivial; hence $\overline{G}=\overline{K}=\overline{N_1}\dots\overline{N_t}\overline{H}$ for $\overline{H}=(C_{\overline{G}}(\overline{N_1})\cap\dots\cap C_{\overline{G}}(\overline{N_t}))$. 
However each $\overline{N_i}\lhd \overline{G}$; so $C_{\overline{G}}(\overline{N_i})$ is normal is $\overline{G}$; which implies $\overline{H}\lhd\overline{G}$. 

Each $\overline{N_i}$ is non-abelian and simple, so $\overline{N}_i\cap C_{\overline{G}}(\overline{N_i})=1$; but as $\overline{H}$ is normal it must contain a minimal normal subgroup, impling that $\overline{H}=1$ and in fact $\overline{G}=\overline{N_1}\dots\overline{N_t}=\overline{N_1}\times\dots\times\overline{N_t}$.

\paragraph{\underline{Step 9: $O_p(G)\not=1$.}}

If $O_p(G)=1$ then $G=N_1\times\dots\times N_t$ by Step 8; as each $N_i$ is a proper normal subgroup of $G$ they must have abelian Sylow $p$-subgroups. Thus the Sylow $p$-subgroups of $G$ must be abelian, contradiction.

\paragraph{\underline{Step 10: $t\leq 2$.}}

If $\overline{G}=\overline{N_1}\times\dots\times \overline{N_t}$ for $t\geq 3$; then by applying Proposition~\ref{SylFromNorm} inductively to $\overline{G}=N_1N_2\dots N_t$ where $\overline{N_i}=N_i/Z(G)$ for $N_i\lhd G$, it follows $G$ must have an abelian Sylow $p$-subgroup.

\subsection{Recovering $G$ as a Central Extension}
Let $M(G)$ denote the Schur multiplier of a group $G$, and $\widetilde{G}$ denote the Schur cover as in \cite{Schur}.

In the previous subsection we established that a minimal non-simple $cl_p$ group $G$ is a central extension in one of the following ways.
\begin{enumerate}
 \item $G/Z(G)\cong S$ for $S$ a simple group with an abelian Sylow $p$-subgroup, or,
 \item $G/Z(G)\cong S_1\times S_2$ for $S_i$ simple groups with abelian Sylow $p$-subgroups.
\end{enumerate}

\subsubsection{Case 1: $G/Z(G)\cong S$}\label{OneSimple}

Throughout this section we will make use of the tables in \cite[Table 8.5]{Schur}, which gives a list of Schur multipliers for all simple groups.

It has been shown that $Z(G)=O_p(G)\not=1$ in Steps 3 and 9, and by \cite[Theorem 33.8]{Asch} $Z(G)$ must be a quotient of $M(S)$.

\paragraph{\underline{Step 1: $p$ must be at most $3$.}}
 
If $p>3$ the only simple groups with Schur multiplier divisible by a prime greater than $3$ are $(A_l(q),C_{(l+1,q-1)})$ and $(^2A_l(q),C_{(l+1,q+1)})$ \cite[Table 8.5]{Schur}, where each pair gives the simple group $S$ and its Schur multiplier $M(S)$.

Recall $A_l(q)$ corresponds to $PSL_{l+1}(q)$ and ${}^2A_l(q)$ to  $PSU_{l+1}(q)$; these simple groups have non trivial Schur multiplier only when $p$ divides $l+1$ and $q-1$ or $q+1$ respectively. 
However in these cases it was established in Section 2 that there is a $p$-element with class size divisible by $p$, giving a contradiction. 

\paragraph{\underline{Step 2: $p$ is not equal to $3$.}}

By \cite[Table 8.5]{Schur} the following families of groups all have 3 dividing the order of $M(S)$.
\begin{center}
\begin{enumerate}
\item $(A_l(q),C_{(l+1,q-1)})$ such that 3 divides $l+1$ and $q-1$
\item $({}^2A_l(q),C_{(l+1,q+1)})$ such that 3 divides $l+1$ and $q+1$
\item $(E_6(q),C_{(3,q-1)})$ and $({}^2E_6(q),C_{(3,q+1)})$. 
\end{enumerate}
\end{center}
In addition to these families; using \cite[Table 8.5]{Schur} and the Atlas \cite{Atlas}, the cases to consider are as follows, where we have only included the simple groups such that $3$ divides the order of $M(S)$ and all $3$-elements have class size coprime to 3.

\begin{enumerate}
 \item $(Alt_6\cong A_1(9),C_6)$,
 \item $(Alt_7,C_6)$,
 \item $(M_{22},C_{12})$,
 \item $(O'N,C_3)$ and
 \item $(A_2(4),C_3\times C_4\times C_4)$.
\end{enumerate}
  
\paragraph{}
We show that each of the above cases leads to a contradiction.

For families $(1)$ and $(2)$, if $3$ divides $n=l+1$ and $q-1$ or $q+1$ for $PSL_n(q)$ and $PSU_n(q)$ respectively; there exists a $3$-element without $3'$-index unless $n=3$ by Section 3.1.
While if $n=3$ and 3 divides $q-1$ or $q+1$ respectively, then $G\cong\widetilde{S}$ is $SL_3(q)$ or $SU_3(q)$ respectively.
In both cases a $3$ element $g=$diag$(\zeta,\zeta^{-1},1)$ can be chosen with class size divisible by $3$.
 
If $S\cong E_6(q)$, then the Weyl group $E_6$ is a subquotient of $E_6(q)$ \cite{SGLT}; and $E_6$ has a non-abelian Sylow 3-subgroup by using the Atlas \cite{Atlas}. 
Thus $E_6(q)$ has a non-abelian Sylow 3-subgroup.
If $S\cong$ $^2E_6(q)$ and $3$ divides $q+1$, by \cite[Sections 4.11 and 4.8.6]{Fsg} there exists the following inclusion, $PSp_{6}(q)\leq F_4(q)\leq {}^2E_6(q)$. By applying the construction of Sylow subgroups in \cite[Section 3]{WeirClass} it can be seen that the Sylow 3-subgroup of $PSp_6(q)$ arises as a non-trivial wreath product so is non-abelian. 
Thus $^2E_6(q)$ has a non-abelian Sylow 3-subgroup. 
 
For the remaining groups, the maximal 3-part of the Schur multiplier is $C_3$; hence if $G/Z(G)$ is one of these simple groups, it follows that $Z(G)\cong C_3$. Thus $G$ must be isomorphic to the triple cover $3.S$ of $S$. 
Using the Atlas \cite{Atlas}, there exists a 3-element $\overline{g}$ in each $S$; such that for $\pi$ the canonical homomorphism from $G$ into $G/Z(G)$ we obtain that $\pi^{-1}(\overline{g}^S)$ in $3.S$ forms a single conjugacy class.
So by Lemma~\ref{ClassSpilt} it follows that $|C_G(g)|_3=|C_S(\overline{g})|_3=|S|_3<|G|_3$ and thus $G$ is not $cl_3$.

\paragraph{}
It has now been established that for a minimal non-simple $cl_p$ group $G$, the quotient $G/Z(G)$ cannot be isomorphic to a simple group $S$.

\subsubsection{Case 2: $G/Z(G)\cong S_1\times S_2$ with $S_i$ simple and non-abelian}
Using \cite[Theorem 2.2.10]{Schur} it can be seen that $M(S_1\times S_2)\cong M(S_1)\times M(S_2)$ and $\widetilde{S_1\times S_2}\cong \widetilde{S_1}\times \widetilde{S_2}$.
Additionally there exists a surjective homomorphism $\alpha:\widetilde{S_1\times S_2}\rightarrow G$, such that $K_\alpha=$ Ker$(\alpha)\leq Z(\widetilde{S_1\times S_2})$ and $1\not=Z(G)=O_p(G)\cong (M(S_1)\times M(S_2))/K_\alpha$ by \cite[Theorem 33.8]{Asch}. 

\paragraph{\underline{Step 1: $p$ must be at most $3$.}}

The prime $p$ divides at least one of $|M(S_1)|$ or $|M(S_2)|$, so without loss of generality assume $|M(S_1)|$ is divisible by $p$.

By the same argument as in Step 1 of Section~\ref{OneSimple}, there exists a $p$-element with class size divisible by $p$ in $S_1$, so the same must hold in $S_1\times S_2$.

\paragraph{\underline{Step 2: $p$ is not equal to $3$.}}
In this case $Z(G)=O_3(G)$ and so $3$ divides the order of $M(S_1)\times M(S_2)$.

By Step 2 of 4.1.1 we can assume that each $S_i$ is not $E_6(q)$, $^2E_6(q)$, $A_l(q)$ or $^2A_l(q)$ such that $l+1>3$, as these simple groups either have a non-abelian Sylow 3-subgroup or there is a $3$-element with 3 dividing its class size.
Thus by looking at the remaining list of possible simple groups after excluding those groups just discussed; then each such remaining group has Schur multiplier not divisible by $9$, hence we obtain two cases for the structure of the 3-part of $M(S_1)\times M(S_2)$, these are $C_3$ and $C_3\times C_3$.

\paragraph{\underline{Case 2.1: The $3$-part of $M(S_1)\times M(S_2)$ is $C_3$.}} Without loss of generality it can be assumed that 3 divides the order of $M(S_1)$ and not the order of $M(S_2)$. In this case $Z(G)\cong C_3$ and it follows that $G\cong 3.S_1\times S_2$. 
However the groups $3.S_1$ are either not $cl_3$ or $S_1$ has a non-abelian Sylow 3-subgroup by Step 2 of Section~\ref{OneSimple}. Thus $S_1\times S_2$ either has a non-abelian Sylow 3-subgroup or $3.S_1\times S_2$ is not a $cl_3$ group.

\paragraph{\underline{Case 2.2: The $3$-part of $M(S_1)\times M(S_2)$ is $C_3\times C_3$.}}
Both $G$ and $3.S_1\times 3.S_2$ are perfect central extensions of $S_1\times S_2$, so there exists homomorphisms $\alpha$ and $\beta$ respectively such that 
\[
 Z(G)\cong M(S_1\times S_2)/K_\alpha \text{ and } Z(3.S_1\times 3.S_2)\cong M(S_1\times S_2)/K_\beta.
\]

Since $M(S_1\times S_2)$ is abelian, then $M(S_1\times S_2)\cong C_3\times C_3\times O_{3'}(M(S_1\times S_2))$.
As $Z(3.S_1\times 3.S_2)$ is a group of order 9, it follows that $K_\beta=O_{3'}(M(S_1\times S_2))$.
The group $Z(G)$ has order either  9 or 3, so $K_\alpha$ must contain $O_{3'}(M(S_1\times S_2))$ and therefore $K_\beta\lhd K_\alpha$.

Appealing to the third isomorphism theorem yields that $G\cong (3.S_1\times 3.S_2)/N$ for $N$ the subgroup in $(3.S_1\times 3.S_2)/K_\beta$ isomorphic to $K_\alpha/K_\beta$, and we see that $N$ has order 1 or 3.

In Step 2 of Section~\ref{OneSimple}, we showed that if the simple group $S_i$ had the order of $M(S_i)$ divisible by 3 and $S_i$ had an abelian Sylow $3$-subgroup; then the group $3.S_i$ contained a $3$-element $g_i$ where $3$ does divide its class size.
Consider the element $\overline{g}$ in $G$ of the form $\overline{g}=(g_1,g_2)N$; then $C_{3.S_1\times 3.S_2}((g_1,g_2))= C_{3.S_1}(g_1)\times C_{3.S_2}(g_2)$, and so $9$ divides $k:=|3.S_1\times 3.S_2|/|C_{3.S_1\times 3.S_2}(g)|$.

At this point we use Proposition~\ref{ClassLift} which tells us that
\[
 |C_G(\overline{g})|=|\{[(g_1,g_2),3.S_1\times 3.S_2]\cap N\}|\cdot |C_{3.S_1\times 3.S_2}((g_1,g_2))|/|N|;
\]
thus the class size of $\overline{g}$ is given by 
\[
 |G|/|C_G(\overline{g})|=k/|\{[(g_1,g_2),3.S_1\times 3.S_2]\cap N\}|.
\]
As $|N|\leq 3$, it follows that $|\{[(g_1,g_2),3.S_1\times 3.S_2]\cap N\}|\leq 3$. Thus $|G|/|C_G(\overline{g})|$ is divisible by 3 and $G$ is not $cl_3$.

\paragraph{}
We have now proven the main theorem which was the aim of this section.
\begin{thm}
Given $G$ a minimal $cl_p$ group then $G$ must be simple.
\end{thm}

Using Corollary~\ref{SubquotientsClp} with Theorem~\ref{SimpleClp}, we have established Theorem~\ref{MainTheorem}.

\bibliographystyle{alpha}
\bibliography{bibfile}

\end{document}